\documentclass[a4paper]{amsart}
\usepackage{enumerate}
\usepackage{pb-diagram}
\usepackage{microtype}
\usepackage{amsmath}
\usepackage{amssymb, pb-diagram}
\usepackage[all]{xy}
\usepackage{graphicx} 
\usepackage{mathabx}
\usepackage{chapterbib}
\usepackage{epsfig}
\usepackage{amsfonts}
\usepackage{amssymb}
\usepackage{amsmath}   
\usepackage{amsthm}
\usepackage{color}
\usepackage{latexsym}
\usepackage{hyperref}
\newtheorem{thm}{Theorem}[section]
\newtheorem*{thmi}{Theorem}

\newtheorem{LM}[thm]{Lemma}
\newtheorem{cor}[thm]{Corollary}
 \theoremstyle{definition}
  \newtheorem{definition}{Definition}[section]

    \newtheorem{rem}[thm]{Remark}

   \DeclareMathOperator {\CAT(0)}{CAT(0)}

    \DeclareMathOperator {\Aut}{Aut}

 \DeclareMathOperator {\rank}{rank} 
 \DeclareMathOperator {\Vol}{Vol}

 \DeclareMathOperator {\Ga}{\Gamma}

  \DeclareMathOperator {\FFF}{F}

     \DeclareMathOperator {\infin}{inf} 
    \DeclareMathOperator {\coh}{H} 
 \DeclareMathOperator {\cd}{cd}

   \DeclareMathOperator {\Z}{\mathbb{Z}}

 \DeclareMathOperator {\FP}{FP} 
 \DeclareMathOperator {\CAT(0)}{CAT(0)}

    \DeclareMathOperator {\CW}{CW}

 \DeclareMathOperator {\Sup}{sup} 
   
 \DeclareMathOperator {\toi}{ \hookrightarrow}    
 \DeclareMathOperator {\tos}{\twoheadrightarrow}

\numberwithin{equation}{section}

\begin{document}

\title[Finiteness properties of non-uniform lattices on $\CAT(0)$ complexes]{Finiteness properties of non-uniform lattices on $\CAT(0)$ polyhedral complexes}

\author{Giovanni Gandini}
\address{School of Mathematics, University of Southampton, Southampton, SO17 1BJ UNITED KINGDOM}


\date{\today}



\begin{abstract}
We show that the homological finiteness length of a non-uniform lattice on a locally finite $\CAT(0)$ $n$-dimensional polyhedral complex is less than $n$. As a corollary, we obtain an upper bound for the homological finiteness length of arithmetic groups over function fields. This gives an easier proof of a result of Bux and Wortman that solved a long-standing conjecture.
\end{abstract}

\maketitle
\section{introduction}
The notion of a lattice in  a locally compact group arises naturally in modern mathematics and  has its  roots in the study of  Lie groups.  A semisimple algebraic group over a local field can be realised as a group of automorphisms of its Bruhat--Tits building, and       their lattices,  called arithmetic lattices are studied since the early 70's. Other examples are given by  tree lattices, which were introduced in the beginning of the  90's by Bass and Lubotzky. Tree lattices are lattices in the isometry group of a locally finite tree \cite{treelattices}.  More recently, lattices in isometry groups of higher dimensional locally finite cell complexes have appeared in literature \cite{Thomas2007, FHT}.

Let $\bold{X}^{n}$ be $S^{n}$, $\mathbb{R}^{n}$ or $\mathbb{H}^{n}$ with Riemannian metrics of constant curvature $1$, $0$ and $-1$ respectively.
A finite-dimensional cell complex $X$ is a polyhedral complex  if each $n$-dimensional open cell  is isometric to the interior of a compact convex polyhedron in $\bold{X}^n$ and the restrictions of the attaching maps are isometries onto its open cells.\
Suppose now that $X$ is a  locally finite $\CAT(0)$ polyhedral complex, then  its group of  cellular isometries  $\Aut(X)$  has a natural structure of a locally compact topological group. We are interested in non-uniform lattices on X. We recall that   a  lattice on $X$ is a discrete subgroup $\Ga$ of $\Aut(X)$ with finite covolume. $\Ga$ is said to be uniform if $\Ga \backslash \Aut(X)$ is compact.  

The homological finiteness length $\phi(G)$  of a group $G$ is a generalisation of the concepts of finite generability and finite presentability. 
The main result of  this paper is a  bound for  the homological finiteness length of certain lattices on $X$.
\begin{thmi} If $\Ga$ is a non-uniform lattice on a locally finite $\CAT(0)$ polyhedral complex of dimension $n$, then   $\phi(\Ga)< n$.
\end{thmi}
The finiteness properties of arithmetic groups have been widely investigated by Abels, Abramenko, Behr, Bux, Serre and Wortman just to mention a few. An upper bound for $\phi(\Ga)$ in the case of  arithmetic groups over function fields is given in \cite{BW}. Our result provides the same bound as a corollary.
\subsection*{Acknowledgements} The author would like to thank Kevin Wortman, Kai--Uwe Bux and Stefan Witzel for the interesting discussions and comments without which this note would not be possible. Moreover, the author is grateful to his supervisor Brita E. A. Nucinkis for her encouragement and advice.

\section{Finiteness properties and $\CAT(0)$ polyhedral complexes}
A cohomological finiteness condition is a group-theoretical property that is satisfied by any group admitting a finite $K(G,1)$.  Since every non-trivial finite group does not admit a  finite-dimensional $K(G,1)$, \emph{being torsion-free} is a cohomological finiteness condition, but not a finiteness condition in the usual group-theoretical sense. 
On the other hand the property of being \emph{locally finite} is a classical but not cohomological  finiteness condition. However there are finiteness conditions that  agree, for example \emph{being finitely generated},   \emph{being finitely presented}... \\
A generalisation of these properties brings us to the concepts of 
cohomological conditions of finite type. More precisely,  a group $\Ga$ is of type $\FP_{n}$  if the trivial $\Z\Ga$-module $\Z$ admits a resolution of finitely generated projective $\Z\Ga$-modules up to dimension $n$. If $\Ga$ is of type $\FP_{n}$ for every $n\geq 0$, then $\Ga$ is said to be of type $\FP_{\infty}$. 
A group is of type $\FFF_{n}$ if it admits a $K(G, 1)$ with finite $n$-skeleton; and $\Ga$ is of type $\FFF_{\infty}$ if it is of type $\FFF_{n}$ for every $n\geq 0$.
For a group being finitely generated is equivalent to being of type $\FP_{1}$. A group is finitely presented if and only if it is of type $\FFF_{2}$.   For $n\geq 2$ a group is of type $\FFF_{n}$ if and only if it is finitely presented and of type $\FP_{n}$.
Bestvina and Brady  show the existence of  non-finitely presented groups of type $\FP_{2}$  \cite{bb-97}.

The \emph{homological finiteness length} of $\Ga$ is defined as $$\phi(\Ga) = \sup\{ m | \,\mbox{$\Ga$ is of type $\FP_{m}$\}}.$$
It is worth mentioning that Abels and Tiemeyer generalise the above finiteness conditions for discrete groups to compactness properties of locally compact groups \cite{abels}.

We begin by recalling the terminology and in doing so we follow closely \cite{Thomas2007} and  \cite{FHT}.
 Let $\bold{X}^{n}$ be $S^{n}$, $\mathbb{R}^{n}$ or $\mathbb{H}^{n}$ with Riemannian metrics of constant curvature $1$, $0$ and $-1$ respectively.
A finite-dimensional $\CW$-complex $X$ is a  \emph{polyhedral complex}  if it satisfies the followings: 
\begin{itemize}
\item each open cell of dimension $n$ is isometric to the interior of a compact convex polyhedron in $\bold{X}^{n}$;
\item
 for each cell $\sigma$ of $X$, the restriction of the attaching map to each open  $\sigma$-face of codimension one  is an isometry onto an open cell of $X$.
\end{itemize}
Let $\Aut(X)$ be the full group of  cellular isometries of $X$.  A subgroup $H\leq \Aut(X)$ acts \emph{admissibly} on $X$ if the  set-wise stabiliser of each cell coincides with its point-wise stabiliser. 
 \begin{rem} Every subgroup $G \leq \Aut(X)$ acts admissibly  on the barycentric subdivision of $X$. Furthermore, if $G$ acts admissibly on a $\CAT(0)$ polyhedral complex, then the fixed point set $X^{G}$ forms a subcomplex of $X$.
 \end{rem}
A subgroup $\Gamma$ of a locally compact topological group $G$ with left-invariant Haar measure $\mu$ is a \emph{lattice} if:
\begin{itemize}
\item $\Gamma$ is discrete, and
\item $\mu(\Gamma \backslash G)< \infty$.
\end{itemize}
 Moreover,  $\Aut(X)$ is locally compact whenever $X$ is locally finite and so it makes sense to talk about lattices on locally finite  $\CAT(0)$ polyhedral complexes.
Let $G$ be a locally compact group with left-invariant Haar measure $\mu$.  Let $\Ga$ be a discrete subgroup of $G$ and $\Delta$  be a $G$-set with  compact and open stabilisers. The $\Delta$-covolume, denoted by $\Vol(\Ga\backslash\backslash \Delta)$, is defined to be $ \sum_{\delta\in \Ga\backslash \Delta} \frac{1}{|\Ga_{\delta}|}\leq \infty$.
\begin{LM}[ \cite{treelattices}, Chapter 1]\label{clas} Let $X$ be a locally finite $\CAT(0)$ polyhedral complex with vertex set $V(X)$. If $\Ga$ is a subgroup of $G=\Aut(X)$, then:
\begin{itemize}
\item $\Gamma$ is discrete if and only if the stabiliser $\Ga_{x}$ is finite for each $x \in V (X)$;
\item $\mu(\Ga\backslash G)<\infty$ if and only if  $\Vol(\Ga\backslash\backslash X)< \infty$.  Moreover,  the Haar measure $\mu$ can be normalised in such  a way that for every discrete $\Ga\leq G$, $\mu(\Ga \backslash G) = \Vol(\Ga\backslash\backslash X$).
\end{itemize}
\end{LM}
\section{Homological finiteness length of non-uniform lattices on $\CAT(0)$ polyhedral complexes}
The next lemma is a well-known criterion for finiteness that follows from Theorems 2.2 and 3.2 in \cite{brownfi}.
\begin{LM}[\cite{brownfi}]\label{brow} Let $\Ga$ be a group that acts on an $n$-dimensional contractible $\CW$-complex with stabilisers of type $\FFF_{\infty}$. Then, $\Ga$ is of type $\FFF_{\infty}$ if and only if it is of type $\FFF_{n}$.
\end{LM}
\begin{definition} The cohomological dimension of $\Ga$ over a ring $R$ is defined as
\begin{align} 
\cd_{R} \Ga &=
\infin\{ n \,|  \;R  \mathrm{ \;admits \;an\;} R\Ga\mbox{-}\mathrm{projective\; resolution \;of\; length} \;n\} \nonumber\\ &=
\Sup\{n \, |\, \coh_{R\Ga}^{n}(G; M)\neq 0,\; \mathrm{for\;some }\,   R\Ga \mbox{-module} \;M \}.  \nonumber
\end{align}
 \end{definition}
The next result, due to Kropholler, is the key ingredient in the proof of Theorem \ref{main}.
\begin{thm}[Proposition, \cite{MR1246274}]\label{krop}
Every group of finite rational cohomological dimension and of type $\FP_{\infty}$  has a bound on the orders of its finite subgroups.
\end{thm}
Furthermore, Kropholler in \cite{MR1246274} shows that every ${\scriptstyle \mathbf H}\mathfrak{F}$-group of type $\FP_{\infty}$  has a bound on the orders of its finite subgroups.
\begin{thm}\label{main} If $\Ga$ is a non-uniform lattice on a locally finite $\CAT(0)$ polyhedral complex of dimension $n$, then   $\phi(\Ga)< n$.
\begin{proof} 
Let $\Ga$ be a non-uniform lattice on a $\CAT(0)$ polyhedral complex $X$ of dimension $n$.   
By Lemma \ref{clas} $\mu(\Ga\backslash\Aut(X))=\sum_{\sigma\in \Ga\backslash X}\frac{1}{|\Ga_{\tilde{\sigma}}|}$,  where  $\sigma= [\tilde{\sigma}]$.
Since $\Ga$ is non-uniform, 
the set $\Ga \backslash X$ is infinite and so  for any $n$ there is some $\sigma \in\Ga \backslash X$ such that $\frac{1}{|\Ga_{\tilde{\sigma}}|}< \frac{1}{n}$. Therefore,  there is no bound on the orders of the stabilisers (which are finite), and so there is no bound on the orders of the finite subgroups of $\Ga$.\\
 In view of Lemma \ref{brow} and Theorem \ref{krop}, it only remains to argue that the rational  cohomological dimension of $\Gamma$ is at most $n$.
 Since every  $\CAT(0)$ space is contractible \cite{BH}, $\Ga$ acts on an $n$-dimensional contractible $\CW$-complex with finite stabilisers.  The augmented cellular chain complex of $X$  is an exact sequence of the form:
$$\bigoplus_{i_{n}\in I_{n}}\Z[\Ga_{i_{n}}\backslash\Ga]\toi \bigoplus_{i_{n-1}\in I_{n-1}}\Z[\Ga_{i_{n-1}}\backslash\Ga]\to \dots \to\bigoplus_{i_{0}\in I_{0}}\Z[\Ga_{i_{0}}\backslash\Ga] \tos\mathbb{Z},$$
where $\Ga_{i_{j}}$ are finite subgroups of $\Ga$ for every $0\leq j \leq n$. Since $\mathbb{Q}$ is flat over $\Z$ and $\mathbb{Q}\otimes \Z [H\backslash \Gamma] \cong \mathbb{Q}[H\backslash \Gamma]$ for any $H\leq \Gamma$, tensoring this sequence with $\mathbb{Q}$ over $\Z$ leads to the exact sequence:
$$\bigoplus_{i_{n}\in I_{n}}\mathbb{Q}[\Ga_{i_{n}}\backslash\Ga]\toi \bigoplus_{i_{n-1}\in I_{n-1}}\mathbb{Q}[\Ga_{i_{n-1}}\backslash\Ga]\to \dots \to\bigoplus_{i_{0}\in I_{0}}\mathbb{Q}[\Ga_{i_{0}}\backslash\Ga] \tos\mathbb{Q}.$$
Now, $\bigoplus_{i_{j}\in I_{j}}\mathbb{Q}[\Ga_{i_{j}}\backslash\Ga]$ is a $\mathbb{Q}\Ga$-projective module  for every $0\leq j \leq n$, and so $\cd_{\mathbb{Q}}\Ga\leq n$.

Theorem \ref{krop} implies that $\Ga$ is not of type $\FP_{\infty}$. Finite groups are of type $\FFF_{\infty}$ and Lemma \ref{brow} completes the proof. Note that this final step can be also achieved by applying Proposition 1 in \cite{bound}.
\end{proof} 
\end{thm}
\begin{rem}Note that if a  finite group acts on a locally finite  $\CAT(0)$ polyhedral complex, then it is contained in the stabiliser of some cell. Now, let $F$ be a finite subgroup of a non-uniform lattice $\Ga$ acting admissibly on a locally finite $\CAT(0)$ polyhedral complex $X$. Since $F$ acts admissibly on $X$, $X^{F}$ is contractible \cite{BH}. In particular, $X$ is a model for $\underline{E}\Ga$. 
\end{rem}
There are not many results that hold for all non-uniform lattices on $\CAT(0)$ polyhedral complexes. As a first immediate application we obtain a classical result.
\begin{cor}
If  $X$ is a tree, then every non-uniform lattice in $\Aut(X)$ is not finitely generated.  
More generally, a   non-uniform lattice on a product of $n$ trees is not of type $\FP_{n}$.
\end{cor}
\begin{cor} Every non-uniform lattice on a locally finite $2$-dimensional $\CAT(0)$ polyhedral complex is not finitely presented.
\end{cor}
Before the last corollary we need to recall some more standard nomenclature. Let $K$ be a global function field, and $S$ be a finite non-empty
set of pairwise inequivalent valuations on $K$.  Let $\mathcal{O}_{S}\leq K$ be the
ring of $S$-integers. Denote a reductive K-group by $\mathbf{G}$. Given a valuation $v$ of $K$, $K_{v}$ is the completion of $K$
with respect to $v$. If $L/K$ is a field extension, the $L$-$\rank$ of $\mathbf{G}$, $\rank_{L} \mathbf{G}$ is the dimension of a maximal $L$-split torus of $\mathbf{G}$. The $K$-group $\mathbf{G}$ is $L$-isotropic if $\rank_{L}\mathbf{G}\neq 0$. As in  \cite{BW}, to any $K$-group $\mathbf{G}$, there is associated a non-negative integer $k(\mathbf{G}, S) = \sum_{v\in S}\rank_{K_{v}} \mathbf{G}$.
We are now ready to state and reprove the Theorem of Bux and Wortman. 


 \begin{cor}[Theorem 1.2, \cite{BW}]\label{1.2}
Let  $\bold{H}$ be a connected non-commutative absolutely almost simple $K$-isotropic $K$-group. Then $\phi(\bold{H}(\mathcal{O}_{S})) \leq k(\bold{H}, S) -1$.
\begin{proof}Let $\bold{H}$  be a connected non-commutative absolutely almost simple $K$-isotropic $K$-group.   Let $H$ be  $\prod_{v\in S}\bold{H}(K_{v})$,  there is  a $k(\bold{H}, S)$-dimensional Euclidean building $X$ associated to $H$. $X$  is a locally finite $\CAT(0)$ polyhedral complex. 
The arithmetic group $\bold{H}(\mathcal{O}_{S})$ becomes  a lattice of $H$ via the diagonal embedding.
$\bold{H}$ is $K$-isotropic if and only if $\bold{H}(\mathcal{O}_{S})$ is non-uniform by \cite{harder}. An application of Theorem \ref{main} completes the proof.
\end{proof}
\end{cor}


\begin{rem} Theorem \ref{main} gives the upper bound  on the homological finiteness length  of arithmetic groups over function fields,  a historical overview can be found in \cite{BW}.
In a recent remarkable paper \cite{bgw} Bux, Gramlich and Witzel   showed that $\phi(\bold{H}(\mathcal{O}_{S})) = k(\bold{H}, S) -1$.  Calculating the homological finiteness length of non-uniform lattices on $\CAT(0)$ polyhedral complexes is an ambitious open problem. We conclude by mentioning that Thomas and Wortman exhibit examples of non-finitely generated non-uniform  lattices on regular right-angled buildings \cite{T-W}. This shows that the upper bound of Theorem \ref{main} is not sharp and in particular, that the Theorem of Bux, Gramlich and Witzel does not hold for all non-uniform lattices on locally finite $\CAT(0)$ polyhedral complexes.
\end{rem}
\bibliographystyle{amsplain}
\bibliography{math}
\end{document}